\documentclass[11pt]{article}
\usepackage{amssymb,latexsym,amsmath}
\topskip  \newtheorem{theorem}{Theorem}[section]

\newtheorem{lemma}[theorem]{Lemma}

\newenvironment{proof}{\medskip\noindent{\it Proof.\ }}{\hfill \mbox{$\Box$}\medskip}

\begin{document}

\pagenumbering{arabic}
\def\llim{\lim_{n\rightarrow\infty}}
\def\ls{\leq}
\def\gs{\geq}
\def\LL{\frak L}
\def\qq{{\bold q}}
\def\txx{{\frac1{2\sqrt{x}}}}
\def\B{\Box}
\def\BB{\Box\hspace{-2.5pt}\Box}
\def\md{\mbox{\,{{\footnotesize mod}}\,}}
\def\th{\mbox{th}}

\title{Squaring the terms of an $\ell^{\th}$ order linear recurrence}
\maketitle

\begin{center}Toufik Mansour \footnote{Research financed by EC's
IHRP Programme, within the Research Training Network "Algebraic
Combinatorics in Europe", grant HPRN-CT-2001-00272}
\end{center}
\begin{center}
Department of Mathematics, Chalmers University of Technology,
S-412~96 G\"oteborg, Sweden

{\tt toufik@math.chalmers.se}
\end{center}
\section{Introduction and the Main result}
An {\em $\ell^{\th}$ order linear recurrence} is a sequence in
which each is a linear combination of the $\ell$ previous terms.
The symbolic representation of an $\ell^{\th}$ order linear
recurrence defined by
\begin{equation}\label{eqflp}
a_{n}=\sum_{j=1}^\ell p_j a_{n-j}=
p_1a_{n-1}+p_2a_{n-2}+\cdots+p_\ell a_{n-\ell},
\end{equation}
is $(a_n(c_0,\ldots,c_{\ell-1};p_1,\ldots,p_\ell))_{n\geq0}$, or
briefly $(a_n)_{n\geq0}$, where the $p_i$ are constant
coefficients, with given $a_j=c_j$ for all $j=0,1,\ldots,\ell-1$,
and $n\geq\ell$; in such a context, $(a_n)_{n\geq0}$ is called an
{\it $\ell$-sequence}.

In the case $\ell=2$, this sequence is called Horadam's sequence
and was introduced, in 1965, by Horadam~\cite{Ho1,Ho2}, and it
generalizes many sequences (see~\cite{HW,HM}). Examples of such
sequences are the Fibonacci numbers $(F_n)_{n\geq0}$, the Lucas
numbers $(L_n)_{n\geq0}$, and the Pell numbers $(P_n)_{n\geq0}$,
when one has the following initial conditions: $p_1=p_1=c_1=1$,
$c_0=0$; $p_1=p_2=c_1=1$, $c_0=2$; and $p_1=2$, $p_2=c_1=1$,
$c_0=0$; respectively. In 1962, Riordan~\cite{R} found the
generating function for powers of Fibonacci numbers. He proved
that the generating function
$\mathcal{F}_k(x)=\sum_{n\geq0}F_n^kx^n$ satisfies the recurrence
relation
$$(1-a_kx+(-1)^kx^2)\mathcal{F}_k(x)=1+kx\sum_{j=1}^{[k/2]}(-1)^j\frac{a_{kj}}{j}\mathcal{F}_{k-2j}((-1)^ix)$$
for $k\geq1$, where $a_1=1$, $a_2=3$, $a_s=a_{s-1}+a_{s-2}$ for
$s\geq3$, and $(1-x-x^2)^{-j}=\sum_{k\geq0} a_{kj}x^{k-2j}$.
Horadam~\cite{Ho2} gave a recurrence relation for
$\mathcal{H}_k(x)$ (see also~\cite{HR}). Haukkanen~\cite{Ha}
studied linear combinations of Horadam's sequences and the
generating function of the ordinary product of two of Horadam's
sequences. Recently, Mansour~\cite{M} found a formula for the
generating functions of powers of Horadam's sequence. In this
paper we interested in studying the generating function for
squaring the terms of the $\ell$-sequence, that is,
$$\mathcal{A}_\ell(x)=\mathcal{A}_\ell(x;c_0,\ldots,c_{\ell-1};p_1,\ldots,p_\ell)
=\sum\limits_{n\geq0}a_n^2(c_0,\ldots,c_{\ell-1};p_1,\ldots,p_\ell)x^n.$$

The main result of this paper can be formulated as follows. Let
$\Delta_\ell=(\Delta_\ell(i,j))_{0\leq i,j\leq \ell-1}$ be the
$\ell\times\ell$ matrix
$$\Delta_\ell(i,j)=\left\{
\begin{array}{ll}
1-\sum_{s=1}^\ell p_j^2x^j,&i=j=0\\
-2xv_j                  ,&i=0\mbox{ and }1\leq j\leq \ell-1\\
-p_{i}x^{i},         &1\leq i\leq\ell-1\mbox{ and }j=0\\
\delta_{i,j}-p_{i-j}x^{i-j}-p_{i+j}x^i, & 1\leq i\leq\ell-1\mbox{
and }1\leq j\leq \ell-i\\
\delta_{i,j},            &1\leq i\leq\ell-1\mbox{ and }
\ell+1-i\leq j\leq \ell-1
\end{array}
\right.$$ where $v_j$ is given by
$$v_j=p_1p_{j+1}+p_2p_{j+2}x+\cdots+p_{\ell-j}p_\ell x^{\ell-j-1},$$
for all $j=1,2,\ldots,\ell-1$, we define $p_i=0$ for $i\leq 0$,
and $\delta_{i,j}=\left\{\begin{array}{ll}
                              1,& \mbox{if } i=j\\
                              0,& \mbox{if } i\neq j
                              \end{array}\right.$
.

Let $\Gamma_\ell=(\Gamma_\ell(i,j))_{0\leq i,j\leq \ell-1}$ be the
$\ell\times\ell$ matrix
$$\Gamma_\ell(i,j)=\left\{
\begin{array}{ll}
x\sum_{s=0}^{\ell-1}(c_s^2-w_{s-1}^2)x^j,   & i=j=0\\
x^{i+1}\sum_{s=0}^{\ell-1-i}c_s(c_{s+i}-w_{s+i-1})x^s,& j=0\mbox{
and }1\leq i\leq\ell-1\\
\Delta_\ell(i,j), & 0\leq i\leq\ell-1\mbox{ and }1\leq j\leq
\ell-1
\end{array}
\right.$$ where $w_j$ is given by
$$w_j=p_1c_j+p_2c_{j-1}+\cdots+p_{j+1}c_0=\sum_{s=1}^{j+1}p_sc_{j+1-s},$$
for $j=0,1,\ldots,\ell-2$ with $w_{-1}=0$.

\begin{theorem}\label{thm}
The generating function $\mathcal{A}_\ell(x)$ is given by
        $$\frac{\det(\Gamma_\ell)}{x\det(\Delta_\ell)}.$$
\end{theorem}

The paper is organized as follows. In Section~\ref{section_2} we
give the proof of Theorem~\ref{thm} and in Section~\ref{section_3}
we give some applications for Theorem~\ref{thm}.
\section{Proofs}\label{section_2}
Let $(a_n)_{n\geq0}$ be a sequence satisfying
Relation~(\ref{eqflp}) and $\ell$ be any positive integer. We
define a family $\{f_d(n)\}_{d=0}^{\ell-1}$ of sequences by
$$f_d(n)=a_{n-1}a_{n-1-s},$$
and a family $\{F_d(x)\}_{d=0}^{\ell-1}$ of generating functions
by
\begin{equation}\label{eqaa}
F_{d}(x)=\sum_{n\geq1}a_{n-1}a_{n-1-d}x^{n}.
\end{equation}
Now we state two relations (Lemma~\ref{lema} and Lemma~\ref{lemb})
between the generating functions $F_{d}(x)$ and
$F_0(x)=x\mathcal{A}_\ell(x)$ that play the crucial roles in the
proof of Theorem~\ref{thm}.

\begin{lemma}\label{lema}
We have
$$F_0(x)=F_0(x)\sum_{j=1}^{\ell}p_j^2x^j+2x\sum_{j=1}^{\ell-1}v_jF_j(x)+x\sum_{j=0}^{\ell-1}(c_j^2-w_{j-1}^2)x^j.$$
\end{lemma}
\begin{proof}
Since the sequence $(a_n)_{n\geq0}$ satisfying
Relation~(\ref{eqflp}) we get that
$$a_n^2=\sum_{j=1}^{\ell}p_j^2a_{n-j}^2+2\sum_{1\leq
i<j\leq\ell}p_ip_ja_{n-i}a_{n-j},$$ for all $n\geq\ell$.
Multiplying by $x^n$ and summing over $n\geq\ell$ together with
the following facts:
\begin{enumerate}
\item
$\sum\limits_{n\geq\ell}a_n^2x^n=\frac{1}{x}\sum\limits_{n\geq\ell}f_0(n+1)x^{n+1}=\frac{1}{x}\left(F_0(x)-\sum\limits_{j=1}^{\ell}a_{j-1}^2x^j\right)$,

\item
$\sum\limits_{n\geq\ell}a_{n-j}^2x^n=\sum\limits_{n\geq\ell}f_0(n-j+1)x^n=x^{j-1}\left(F_0(x)-\sum\limits_{t=1}^{\ell-j}a_{t-1}^2x^t\right)$,

\item
$\sum\limits_{n\geq\ell}a_{n-i}a_{n-j}x^n=\frac{1}{x}\sum\limits_{n\geq\ell+1}f_{j-i}(n-i)x^n=x^{i-1}\left(F_{j-i}(x)-\sum\limits_{d=j-i+1}^{\ell-i}a_{d-1}a_{d-j+i-1}x^d\right)$,
\end{enumerate}
we have that
$$\begin{array}{ll}
F_0(x)&=F_0(x)\sum\limits_{j=1}^\ell p_j^2x^j+2\sum\limits_{1\leq
i<j\leq\ell}p_ip_j x^iF_{j-i}(x)\\
&\qquad+\sum\limits_{j=1}^\ell a_{j-1}^2x^j-
\sum\limits_{j=1}^{\ell}\sum\limits_{i=1}^{\ell-j}p_j^2a_{i-1}^2x^{j+i}-
2\sum\limits_{1\leq i<j\leq\ell}\sum\limits_{d=j-i+1}^{\ell-i}
p_ip_ja_{d-1}a_{d-(j-i)-1}x^{i+d}\\
&=F_0(x)\sum\limits_{j=1}^\ell
p_j^2x^j+2x\sum_{j=1}^{\ell-1}v_jF_j(x)+x\sum\limits_{j=0}^{\ell-1}(a_j^2-w_{j-1}^2)x^j.
\end{array}$$
Hence,using the fact that $a_j=c_j$ for $j=0,1,\ldots,\ell-1$ we
obtain the desired result.
\end{proof}

\begin{lemma}\label{lemb}
For any $i=1,2,\ldots,\ell-1$,
$$F_i(x)=p_ix^iF_0(x)+\sum_{j=1}^{\ell-i}(p_{i-j}x^{i-j}+p_{i+j}x^{i})F_j(x)+x^{i+1}\sum_{j=0}^{\ell-1-i}c_j(c_{i+j}-w_{i+j-1})x^j.$$
\end{lemma}
\begin{proof}
By direct calculations we have for $n\geq\ell+1$,
$$f_i(n)=a_{n-1}a_{n-1-i}=\sum_{j=1}^{\ell}p_ja_{n-1-j}a_{n-1-i};$$
equivalently,
$$\begin{array}{l}
f_i(n)\\
\,=p_1f_{i-1}(n-1)+p_2f_{i-2}(n-2)+\cdots+p_if_0(n-i)+p_{i+1}f_1(n-i)+\cdots+p_\ell
f_{\ell-i}(n-i). \end{array}$$ As in Lemma~\ref{lema}, multiplying
by $x^n$ and summing over $n\geq\ell+1$ we get
$$\begin{array}{l}
F_i(x)-\sum\limits_{j=i+1}^\ell
a_{j-1}a_{j-1-i}x^j=\sum\limits_{j=1}^i
p_jx^j\left(F_{i-j}(x)-\sum\limits_{d=i-j+1}^{\ell-j}
a_{d-1}a_{d-(i-j)-1}x^d\right)\\
\qquad\qquad\qquad\qquad\qquad\qquad\qquad+\sum\limits_{j=i+1}^\ell
p_jx^i\left(F_{j-i}(x)-\sum\limits_{d=j-i+1}^{\ell-i}a_{d-1}a_{d-(j-i)-1}x^d\right)
\end{array}$$
The rest is easy to check from the definitions.
\end{proof}

\begin{proof}(Theorem~\ref{thm})
Using the above lemmas together with the definitions we have
$$\Delta_k\cdot [F_0(x), F_1(x),F_2(x),\ldots,F_{\ell-1}(x)]^{\mbox{T}}={\bf w}^{\mbox{T}},$$
where the vector ${\bf w}$ is given by
$$\left[\begin{array}{c}
x\sum_{j=0}^{\ell-1}(c_j^2-w_{j-1})x^j\\
x^2\sum_{j=0}^{\ell-2}c_j(c_{j+1}-w_j)x^j\\
x^3\sum_{j=0}^{\ell-3}c_j(c_{j+2}-w_{j+1})x^j\\
\vdots\\
x^{\ell-1}\sum_{j=0}^0c_j(c_{j+\ell-1}-w_{j+\ell-2})x^j
\end{array}\right].$$ Hence, the solution of the above equation gives the
generating function
$F_0(x)=\frac{\det(\Gamma_\ell)}{\det(\Delta_\ell)}$;
equivalently,
$\mathcal{A}_\ell(x)=\frac{\det(\Gamma_\ell)}{x\det(\Delta_\ell)}$,
as claimed in Theorem~\ref{thm}.
\end{proof}
\section{Applications}\label{section_3}
In this section we present some applications of Theorem~\ref{thm}.

{\bf Fibonacci numbers}. Let $F_{k,n}$ be the $n^{\th}$
$k$-Fibonacci number which is given by
$$F_{k,n}=\sum_{j=1}^k F_{k,n-j},$$
for $n\geq k$, with $F_{k,0}=0$ and $F_{k,j}=1$ for
$j=1,2,\ldots,k-1$; in such a context, $F_{2,n}$, $F_{3,n}$, and
$F_{4,n}$ are usually called the $n^{\th}$ Fibonacci numbers,
tribonacci numbers, and tetranacci numbers; respectively. Using
Theorem~\ref{thm} with $c_0=0$ and
$$c_1=c_2=\cdots=c_{k-1}=p_1=p_2=\cdots=p_k=1$$ gives the generating
function $\sum_{n\geq0}F_{k,n}^2x^n$ (see Table~1).

\begin{table}[h]
\begin{center}
\begin{tabular}{|l|l|l|} \hline
  $k$ & The generating function $\sum_{n\geq0}F_{k,n}^2x^n$
  \\ \hline\hline
  $2$ & $\frac{x(1-x)}{(1+x)(1-3x+x^2)}$ \\[5pt]
  $3$ & $\frac{x(1-x-x^2-x^3)}{(1+x+x^2-x^3)(1-3x-x^2-x^3)}$ \\[5pt]
  $4$ & $\frac{x(1-x-5x^2-2x^3-x^4-2x^5+3x^7+x^8)}{1-2x-4x^2-5x^3-8x^4+4x^5+6x^6+x^8-x^{10}}$\\[5pt]
  $5$ & $\frac{x(1-x-5x^2-12x^3-8x^4-10x^5-7x^6-17x^7-8x^8+13x^9+10x^{10}+3x^{11}+9x^{12}+4x^{13})}{1-2x-4x^2-7x^3-11x^4-16x^5+4x^6+7x^7+4x^8+4x^9+7x^{10}-x^{12}-x^{13}-x^{15}}$\\[3pt] \hline
\end{tabular}
\caption{The generating function for the square of the
$k^{\th}$-Fibonacci numbers}
\end{center}
\end{table}

From Table~1, for $k=3$ we obtain
$$\sum_{n\geq0}nF_{3,n}^2x^n=
\frac{x(1-2x+2x^2+12x^3+8x^5+2x^6+4x^7+3x^8+2x^9)}{(x^3-x^2-x-1)^2(x^3+x^2+3x-1)^2}.$$

{\bf Pell numbers}. Let $P_{k,n}$ be the $n^{\th}$ $k$-Pell number
which is given by
$$P_{k,n}=2P_{k,n-1}+\sum_{j=2}^k P_{k,n-j},$$
for $n\geq k$, with $P_{k,j}=1$ for $j=0,1,\ldots,k-1$; in such a
context, $P_{2,n}$ is usually called the $n^{\th}$ Pell number.
Using Theorem~\ref{thm} with $c_j=1$ for $j=0,1,\ldots,k-1$ and
$p_j=1$ for $j=1,2,\ldots,k$ gives the generating function
$\sum_{n\geq0}P_{k,n}^2x^n$ (see Table~2).

\begin{table}[h]
\begin{center}
\begin{tabular}{|l|l|l|} \hline
  $k$ & The generating function $\sum_{n\geq0}P_{k,n}^2x^n$
  \\ \hline\hline
  $2$ & $\frac{1-4x-x^2}{(1+x)(1-6x+x^2)}$ \\[6pt]
  $3$ & $\frac{1-4x-11x^2-13x^3-5x^4-4x^5}{(1-6x-3x^2-x^3)(1-x+2x^2-x^3)}$ \\[6pt]
  $4$ & $\frac{1-4x-12x^2-25x^3-29x^4-3x^5-9x^6-12x^7+13x^8+9x^9}{(1-5x-8x^2-13x^3-20x^4+2x^5+14x^6+x^7+x^8-x^{10})}$\\[6pt]
  $5$ & $\frac{(1+x)(9x^{13}+4x^{12}+2x^{11}+13x^{10}+6x^9-26x^8-6x^7-2x^6-5x^5-14x^4-9x^3-3x^2-2x+1)}{1-2x-4x^2-7x^3-11x^4-16x^5+4x^6+7x^7+4x^8+4x^9+7x^{10}-x^{12}-x^{13}-x^{15}}$\\[3pt] \hline
\end{tabular}
\caption{The generating function for the square of the
$k^{\th}$-Pell numbers}
\end{center}
\end{table}

From Table~2, for $k=2$ we have
$$\sum_{n\geq0}nP_{2,n}^2x^n=
\frac{x(1-2x+10x^2-2x^3+x^4)}{(x+1)^2(x^2-6x+1)^2}.$$

\noindent {\sc 2000 Mathematics Subject Classification}: Primary
11B39; Secondary 05A15
\end{document}